%% file: main.tex
\documentclass[Journal,letterpaper,InsideFigs,SingleSpace,NoLineNumbers]{ascelike-new}
\usepackage[colorlinks,linkcolor=blue,citecolor=blue,urlcolor=blue,draft=false]{hyperref}
\usepackage[utf8]{inputenc}
\usepackage[T1]{fontenc}
\usepackage{lmodern}
\usepackage{graphicx}
\usepackage[figurename=Fig.,labelfont=bf,labelsep=period]{caption}
\usepackage{subcaption}
\usepackage{amsmath}
\usepackage{amssymb}
\usepackage{bm}
\usepackage{doi}
\usepackage{mathtools}
\usepackage{newtxtext,newtxmath}
\usepackage[round]{natbib}
\usepackage{siunitx}
\usepackage{xcolor}
\usepackage{fancyvrb}
\usepackage{cprotect}

\NameTag{Shaw, \today}

\DeclareMathOperator{\E}{\mathbb{E}}
\DeclareMathOperator{\sech}{sech}

\newcommand{\dee}{\mathrm{d}}

\newcommand{\diff}{\:\dee}
\newcommand{\Ensemble}[1]{\left\langle #1 \right\rangle}
\newcommand{\etamean}{\mean{\eta}}
\newcommand{\flow}{\vect{U}}

\newcommand{\flowKmodified}{\flow^{K,\star}}
\newcommand{\flux}{\vect{F}}

\newcommand{\hKmodified}{h^{K,\star}}
\newcommand{\hump}{r}
\newcommand{\humpmean}{\mean{r}}
\newcommand{\Mag}[1]{\left\lvert #1 \right\rvert}
\newcommand{\mean}[1]{\overline{#1}}
\newcommand{\palt}{s}
\newcommand{\pcbasis}{\Phi}
\newcommand{\pcbasisvect}{\bm{\Phi}}
\newcommand{\qKmodified}{q^{K,\star}}
\newcommand{\randomroot}{\hat{\xi}}
\newcommand{\riemannflux}{\vect{\tilde{F}}}

\newcommand{\source}{\vect{S}}
\newcommand{\T}{\intercal}
\newcommand{\vect}{\mathbf}
\newcommand{\velocity}{v}
\newcommand{\zmean}{\mean{z}}
\newcommand{\zmodified}{z^\star}
\newcommand{\zrectbump}{\breve{z}}

\newcommand{\rev}[1]{#1}

\author[1]{James Shaw, Ph.D}
\author[2]{Georges Kesserwani, Ph.D}

\affil[1]{Research Associate, Department of Civil and Structural Engineering, The University of Sheffield, Western Bank, Sheffield S10 2TN, U.K. Email: js102@zepler.net}
\affil[2]{Research Fellow, Department of Civil and Structural Engineering, The University of Sheffield, Western Bank, Sheffield S10 2TN, U.K. Email: g.kesserwani@sheffield.ac.uk}

\title{\rev{Stochastic Galerkin finite volume shallow flow model: well-balanced treatment over uncertain topography}}

\begin{document}

\maketitle

\input{abstract}
\input{intro}
\input{method}
\input{cproperty}
\input{experiments}
\input{conclusions}

\input{ack}
\appendix
\input{software}
\bibliography{references}

\end{document}

%% file: abstract.tex
\begin{abstract}
Stochastic Galerkin methods can quantify uncertainty at a fraction of the computational expense of conventional Monte Carlo techniques, but such methods have rarely been studied for modelling shallow water flows.
Existing stochastic shallow flow models are not well-balanced and their assessment has been limited to stochastic flows with smooth probability distributions.
This paper addresses these limitations by formulating a one-dimensional stochastic Galerkin shallow flow model using a low-order Wiener-Hermite Polynomial Chaos expansion with a finite volume Godunov-type approach, incorporating the surface gradient method to guarantee well-balancing.
Preservation of a lake-at-rest over uncertain topography is verified analytically and numerically.
The model is also assessed using flows with discontinuous and highly non-Gaussian probability distributions.
Prescribing constant inflow over uncertain topography, the model converges on a steady-state flow that is subcritical or transcritical depending on the topography elevation.
Using only four Wiener-Hermite basis functions, the model produces probability distributions comparable to those from a Monte Carlo reference simulation with 2000 iterations, while executing about 100 times faster.
Accompanying model software and simulation data is openly available online.
\end{abstract}

%This paper addresses these limitations by formulating a well-balanced one-dimensional stochastic Galerkin shallow flow model using a low-order Wiener-Hermite \rev{Polynomial Chaos} basis.
%The model is based on a stochastic reformulation of a finite volume Godunov-type approach, incorporating the surface gradient method to guarantee well-balancing.
%Despite its relatively low-order \rev{Wiener-Hermite} basis, the model produces probability distributions comparable to those from a Monte Carlo reference simulation, and executes about 100 times faster.

%% file: intro.tex
\section{Introduction}
Shallow water flows can be highly sensitive to uncertainties in input data: \rev{uncertainties in rating curve extrapolation}, the choice of friction coefficient, \rev{digital elevation model errors and river survey measurement errors} can substantially alter simulated distributions of flood extent, water depth and flow velocity \citep{bates2014,kim-sanders2016,jung-merwade2012,montanari-dibaldassarre2013}.
Conventional uncertainty quantification methods, including the Generalized Likelihood Uncertainty Estimation (GLUE) method \citep{beven-binley1992}, are based on Monte Carlo sampling with randomised inputs, making conventional methods computationally expensive and slow to converge on flow statistics.
As a result, Monte Carlo simulations of flood events are severely constrained by available computing resources, and many sources of uncertainty must be neglected to make probabilistic simulations feasible \citep{neal2013}.
Stochastic Galerkin methods offer a computationally efficient alternative: repeated sampling is eliminated, and convergence is typically far more rapid than Monte Carlo methods \citep{xiu2009,ge2008}.
%Despite widespread interest in stochastic Galerkin and other polynomial chaos methods, existing stochastic Galerkin shallow water models have lacked any well-balanced treatment, and have only been validated for flows with smooth probability distributions \citep{ge2008}.
%This paper addresses these limitations by formulating a well-balanced stochastic Galerkin shallow water model, and assessing the model using flows with discontinuous, non-Gaussian probability distributions.

Polynomial Chaos methods work by introducing to the governing deterministic equations an additional, stochastic dimension, and approximating it with a \rev{so-called} Polynomial Chaos expansion.
Polynomial Chaos methods can be classified as non-intrusive or intrusive.
Non-intrusive methods repeatedly sample a deterministic model with different input values, using the outputs to construct a stochastic solution.
While non-intrusive polynomial chaos methods involve repeated sampling, they require far fewer samples than Monte Carlo methods \citep{ge2008,ge2011}.
Intrusive methods reformulate a deterministic model to produce a stochastic formulation that must be implemented as a new stochastic solver.
The stochastic Galerkin method is an intrusive method which makes a Galerkin projection in stochastic space to produce a decoupled system of equations that are solved in a single model run \rev{to directly evolve the coefficients in the Polynomial Chaos expansion}.
Conventional stochastic Galerkin methods span the entire stochastic space using a single element with a basis chosen from the Askey scheme of orthogonal polynomials \citep{xiu-karniadakis2002}.
When the stochastic flow is sufficiently smooth and the probability distributions are well-represented by the Polynomial Chaos basis, then Polynomial Chaos methods can achieve exponential convergence in the stochastic dimension \citep{xiu-karniadakis2003}.
\citet{sattar-elbeltagy2017} applied a conventional stochastic Galerkin method to the one-dimensional water hammer equations, and their model was able to converge on solutions with smooth, Gaussian probability distributions using a first-order Polynomial Chaos basis with just two basis functions.

Conventional stochastic Galerkin methods have two main shortcomings.
First, strongly nonlinear flows produce steep gradients or discontinuities in stochastic space that are often poorly represented by global Polynomial Chaos bases using a single element \citep{pettersson2014}.
Second, even when exponential convergence is achieved early in the simulation, solution accuracy can deteriorate severely over long simulations \citep{gerritsma2010}.
These shortcomings have motivated a wealth of research into better alternatives, including time-dependent basis functions \citep{gerritsma2010}, multiresolution wavelet decompositions \citep{lemaitre2004a}, multi-element discretisations and adaptive meshing of stochastic space \citep{wan-karniadakis2006,tryoen2010a,pettersson2014,li-stinis2015}.
%Similar adaptive methods have also been developed for non-intrusive stochastic collocation methods, and these are reviewed by \citet{bhaduri2018}.

While stochastic Galerkin methods have received a great deal of attention, few have applied the method to the shallow water equations.
\citet{ge2008} made the first step in this direction, formulating a one-dimensional stochastic Galerkin shallow water model to study solitary wave propagation over uncertain topography.
Their stochastic Galerkin model was 50 times more efficient than Monte Carlo with similar accuracy, demonstrating that stochastic Galerkin shallow water models offer a viable alternative to Monte Carlo methods.
However, their \rev{numerical tests} were limited to flows with smooth, unimodal probability distributions that are more easily represented by stochastic Galerkin methods.
Furthermore, their stochastic Galerkin formulation used a centred difference approximation of the bed slope source term that is not well-balanced.

A well-balanced stochastic Galerkin model was formulated by \citet{jin2016} and applied to scalar nonlinear equations including the inviscid Burgers' equation.
The model inherited the well-balancing property from the underlying deterministic model that used the `interface method' by \citet{jin2001}.
\citet{jin2016} note that the intrusive stochastic Galerkin method has a particular advantage over non-intrusive stochastic collocation:
a stochastic collocation method can only guarantee well-balancing at the collocation nodes where the deterministic well-balanced model is sampled.
In contrast, the stochastic Galerkin model guarantees well-balancing across the entire stochastic space, irrespective of the order of Polynomial Chaos basis.

When applied to the shallow water equations, the well-balanced interface method is equivalent to the surface gradient method by \citet{zhou2001}.
A stochastic reformulation of the surface gradient method is straightforward because it is a linear method involving only linear combinations of the discrete flow variables and topography.
Popular alternatives to the surface gradient method include well-balanced approaches that are also depth-positivity preserving \citep{audusse2004,liang-marche2009}.
Such approaches use nonlinear $\max(\cdot, \cdot)$ operators which greatly complicate their stochastic Galerkin reformulation for realistic applications.

This paper presents a well-balanced stochastic Galerkin model of the one-dimensional shallow water equations \rev{for simulating probabilistic flows that account for measurement errors in flood plain topography or open channel bathymetry.}
\rev{This one-dimensional stochastic model represents the first step towards a two-dimensional probabilistic flood model, and it is also relevant to one-dimensional open channel flows.}
\rev{Being based on a Godunov-type shock-capturing method, the model is capable of simulating subcritical, supercritical and transcritical flows.}
\rev{The model uses a stochastic reformulation of the surface gradient method \citep{zhou2001}, and the well-balancedness of this reformulation} is studied theoretically and verified numerically for \rev{an idealised, motionless lake-at-rest over an uncertain and irregular bed}.
Next, a steady-state test is proposed that provokes a strongly nonlinear flow response over an uncertain bed, resulting in probability distributions that are highly non-Gaussian and discontinuous.
Probability distributions produced by the stochastic Galerkin model are validated against a Monte Carlo reference simulation, and the relative computational performance is discussed.
\rev{Finally, the stochastic Galerkin model is verified for a more realistic steady-state flow over a highly irregular bed with uncertainties that characterise measurement errors in bed elevation.}
Numerical simulation data \citep{shaw-kesserwani2019a} and a Python 3 implementation of the stochastic Galerkin shallow flow model \citep{shaw-kesserwani2019b} are available to download from Zenodo.
Instructions for running the models and interpreting the data are provided in the appendix.

%This paper makes the next step by formulating a well-balanced stochastic Galerkin shallow water model, and assessing the model using flows with discontinuous, non-Gaussian probability distributions.

%% file: method.tex
\section{Deterministic and stochastic shallow flow models}

In this section, a certain deterministic numerical solver of the one dimensional (1D) shallow water equations is outlined in the framework of a finite volume Godunov-type method.
The selected deterministic solver relies on the surface gradient method \citep{zhou2001} to both ensure a well-balanced topography integration and extendibility of the well-balanced property into the stochastic Galerkin case. 
Accordingly, a stochastic Galerkin reformulation is devised that is theoretically well-balanced with uncertain topography under a lake-at-rest hypothesis.

The mathematical model of the shallow water equations represent mass and momentum
conservation principals, and is used in the following conservative form when solving it within a finite volume Godunov-type framework \citep{toro-garcianavarro2007}:
\begin{align}
\frac{\partial \flow(x, t)}{\partial t} + \frac{\partial \flux(\flow(x, t))}{\partial x} = \source(\flow(x, t), z(x)) \label{eqn:swe}
\end{align}
where $\flow = \left[ h, q \right]^\T$ is the flow vector including water \rev{depth} $h$ ($\mathrm{L}$) and unit-width discharge $q = h\velocity$ ($\mathrm{L}^2/\mathrm{T}$) in which $\velocity$ represents the depth-averaged velocity ($\mathrm{L}/\mathrm{T}$), $\flux = \left[ q,  q^2/h + gh^2/2 \right]^\T$ is the flux vector in which $g$ represents the gravitational constant and $\source = \left[ 0, -gh \: \dee z / \dee x \right]^\T$ is the source term vector in which the gradient of the topography $z(x)$ is involved.
\rev{Equation~\eqref{eqn:swe} represents hydraulic flow in an idealised, frictionless channel with a rectangular cross-section of unit width.}

\subsection{Deterministic model}

On a uniform 1D mesh with $M$ elements each of size $\Delta x$, the first-order finite volume method leads to the following discrete element-wise formulation of \rev{the shallow water equations given by equation~\eqref{eqn:swe}}:
\begin{align}
    \flow_i^{(n+1)} = \flow_i^{(n)} - \Delta t
    \left(
    \frac{\riemannflux_{i+1/2}^{(n)} - \riemannflux_{i-1/2}^{(n)}}{\Delta x}
    - \source_i^{(n)} \right) \label{eqn:swe-discrete}
\end{align}
in which $\flow_i^{(n)} = \left[ h_i^{(n)}, q_i^{(n)} \right]^\T$ is a piecewise-constant discretisation of the flow vector at element $i$ and time level $(n)$, and $\riemannflux_{i+1/2}^{(n)}$ is a numerical flux function for linking nonlinear discontinuities associated with the flow vector data at interface $i+1/2$ located between element $i$ and $i + 1$.
Namely, $\riemannflux_{i+1/2}^{(n)} = \riemannflux(\flow_{i+1/2}^-, \flow_{i+1/2}^+)$ where $\flow_{i+1/2}^-$ is the limit of the solution from the side of element $i$, and $\flow_{i+1/2}^+$ is the limit of the solution from the side of element $i+1$.
Within the scope of this work involving a first-order accurate solver, these limits become $\flow_{i+1/2}^- = \flow_i$ and $\flow_{i+1/2}^+ = \flow_{i+1}$ that are used in a numerical flux function based on the Roe approximate Riemann solver \citep{roe-pike1984}.
\rev{Consequently, the deterministic model is able to capture the occurrence of shocks, and simulate subcritical, supercritical and transcritical flows.}

The surface gradient method essentially reconstructs an averaged topography at interface $i+1/2$ that is shared by both elements $i$ and $i+1$, as $\zmodified_{i+1/2} = (z_i + z_{i+1})/2$.
From the reconstructed topography $\zmodified_{i+1/2}$, consistent flow variable limits are accordingly reconstructed based on the actual free-surface elevation data, i.e. $\eta_{i+1/2}^- = h_i^{(n)} + z_i$ and $\eta_{i+1/2}^+ = h_{i+1}^{(n)} + z_{i+1}$, and velocity data, i.e. $\velocity_{i+1/2}^- = q_i^{(n)}/h_i^{(n)}$ and $\velocity_{i+1/2}^+ = q_{i+1}^{(n)} / h_{i+1}^{(n)}$, as: $\hKmodified_{i+1/2} = \eta_{i+1/2}^K - \zmodified_{i+1/2}$ and $\qKmodified_{i+1/2} = \hKmodified_{i+1/2} \velocity_{i+1/2}^K$ (where $K = + \text{ or } -$).
For clarity of presentation, the time level denoted by superscript $(n)$ is omitted from all reconstructed variables.
These reconstructions form new Riemann states $\flowKmodified_{i+1/2} = \left[ \hKmodified_{i+1/2}, \qKmodified_{i+1/2} \right]^\T$ for use to evaluate $\riemannflux_{i+1/2}^{(n)}$.
By analogy, new Riemann limits $\flowKmodified_{i-1/2} = \left[ \hKmodified_{i-1/2}, \qKmodified_{i-1/2} \right]^\T$ at $i - 1/2$ are produced for use to evaluate $\riemannflux_{i-1/2}^{(n)}$.
From the reconstructed limits, a well-balanced discretisation of the source term vector can be produced:
\begin{align}
	\source_i^{(n)} = \left[ 0, -g
	\left( \frac{h^{+,\star}_{i-1/2} + h^{-,\star}_{i+1/2}}{2} \right)
	\left( \frac{z^\star_{i+1/2} - z^\star_{i-1/2}}{\Delta x} \right)
	\right]^\T
	\label{eqn:source}
\end{align}
The well-balanced deterministic model presented in equations~\eqref{eqn:swe-discrete} and \eqref{eqn:source} is used for Monte Carlo simulations, and the deterministic model is also the starting point for a stochastic Galerkin reformulation.

\subsection{Fundamental properties of the Polynomial Chaos basis}

Before presenting the stochastic Galerkin reformulation, it is necessary to consider a single random variable \rev{$A(\theta)$} that maps from the random event $\theta$ to an arbitrary probability distribution with finite variance.
This random variable can be approximated by a Wiener-Hermite Polynomial Chaos expansion \citep{xiu-karniadakis2002}.
The expansion is based on a standard Gaussian random variable $\xi(\theta) \in [-\infty, +\infty]$ having zero mean and unit variance.
The random variable of interest, \rev{$A(\theta)$}, is then approximated as
\rev{\begin{align}
A(\theta) \approx \sum_{p=0}^P A_p \pcbasis_p(\xi(\theta))
\end{align}}
where \rev{$\vect{A} = \left[ A_0, \ldots, A_P \right]^\T$} are the expansion coefficients and $\pcbasisvect = \left[ \pcbasis_0, \ldots, \pcbasis_P \right]^\T$ is the probabilists' Hermite polynomial basis having basis function $\pcbasis_p$ of degree $p$,
\begin{align}
    \pcbasis_p(\xi) = \left( -1 \right)^p \exp \left(\frac{\xi^2}{2}\right)
    \frac{\dee^p}{\dee \xi^p} \exp \left(- \frac{\xi^2}{2} \right)
\end{align}
where $\pcbasis_0 = 1, \pcbasis_1 = \xi, \pcbasis_2 = \xi^2 - 1, \pcbasis_3 = \xi^3 - 3\xi$ and so on.
As the basis order $P$ is increased, the Wiener-Hermite Polynomial Chaos approximation converges on the true random variable \rev{$A(\theta)$} \citep{xiu-karniadakis2002}.

\subsubsection*{Basis orthogonality and commutativity}
The Wiener-Hermite basis $\pcbasisvect$ is orthogonal such that
\begin{align}
	\Ensemble{\pcbasis_p \pcbasis_s} = \Ensemble{\pcbasis_p^2} \delta_{ps}
\end{align}
where $\Ensemble{\cdot}$ is the ensemble average operator and $\delta_{ps}$ is the Kronecker delta that is equal to one when $p = s$ and zero otherwise.
The ensemble average operator is defined as the weighted integral over the standard Gaussian random variable $\xi$:
\begin{align}
	\Ensemble{\alpha(\xi)} = \int_{-\infty}^\infty \alpha(\xi) W(\xi) \diff \xi \label{eqn:ensemble-average}
\end{align}
where $\alpha(\xi)$ is an expression involving any combination of random variables or basis functions, and the weighting function $W(\xi)$ is the standard Gaussian probability density function
\begin{align}
	W(\xi) = \frac{1}{\sqrt{2\pi}} \exp \left(-\frac{\xi^2}{2}\right)
\end{align}
This weighting function ensures that, when $\alpha$ is independent of $\xi$, the ensemble average $\Ensemble{\alpha} = \alpha$.
Finally, the ensemble average of a product of basis functions is commutative such that
\begin{align}
    \Ensemble{\pcbasis_p \pcbasis_s} = \Ensemble{\pcbasis_s \pcbasis_p}
    \label{eqn:commutative}
\end{align}
The commutative property is needed later when verifying the well-balanced property with uncertain topography.

\subsubsection*{Stochastic Galerkin projection of a random variable}
Given the random variable \rev{$A(\theta) = \sum_{p=0}^P A_p \pcbasis_p(\xi(\theta))$}, its Galerkin projection onto a basis function $\pcbasis_l$ with $l = 0, \ldots, P$ is achieved using the ensemble average operator such that, due to orthogonality,
\rev{\begin{align}
	\Ensemble{A(\theta) \pcbasis_l} = A_l \Ensemble{\pcbasis_l^2} \label{eqn:orthogonal}
\end{align}}
where \rev{$A_l$} is the $l$\textsuperscript{th} order expansion coefficient.
Also note that the Galerkin projection of a basis function $\pcbasis_l$ and two random variables, \rev{$A(\theta)$ and $B(\theta)$}, is distributive:
\rev{\begin{align}
	\Ensemble{\left(A(\theta) + B(\theta)\right) \pcbasis_l}
	=
	\Ensemble{A(\theta) \pcbasis_l} + \Ensemble{B(\theta) \pcbasis_l} \label{eqn:distributive}
\end{align}}

\subsubsection*{Mean, variance and high-order moments}
The mean, variance and high-order moments can be calculated for \rev{$A(\theta)$}.
The $m$\textsuperscript{th} moment \rev{$\mu_m[A]$} is defined as
\rev{\begin{align}
\mu_m[A] = \int_{-\infty}^\infty \left(A - \beta \right)^m W(\xi) \diff \xi
    =
    \Ensemble{\left( A - \beta \right)^m} \label{eqn:moment}
\end{align}}
where $\beta = 0$ when $m = 1$ and \rev{$\beta = \mu_1[A]$} for higher-order moments.
Therefore, the mean \rev{$\mu_1[A] = \Ensemble{A} = \sum_{p=0}^P A_p \Ensemble{\pcbasis_p}$}.
Since $\Ensemble{\pcbasis_0} = 1$ and $\Ensemble{\pcbasis_p} = 0$ for $p > 0$ then
\rev{\begin{align}
\mu_1[A] = A_0
\label{eqn:mean}
\end{align}}
The shorthand notation for the mean of \rev{$A$} is \rev{$\mean{A}$}, also known as the expected value, \rev{$\E\left[A\right]$}.

The variance \rev{$\mu_2[X]$} can be derived using the fact that \rev{$\E\left[ \left( A - \E[A] \right)^2 \right] = \E[A^2] - \E^2[A]$}, hence \rev{$\mu_2[A] = \left(\sum_{p=0}^P A_p^2 \Ensemble{\pcbasis_p^2}\right) - A_0^2 \Ensemble{\pcbasis_0}^2$}.
Since $\Ensemble{\pcbasis_0}^2 = \Ensemble{\pcbasis_0^2}$ then
\rev{\begin{align}
    \mu_2[A] &= \sum_{p=1}^P A^2_p \Ensemble{\pcbasis_p^2} \label{eqn:variance}
\end{align}}
The shorthand notation for the variance of \rev{$A$} is \rev{$\sigma^2_A$} and the standard deviation of \rev{$A$} is \rev{$\sigma_A$}.

\subsubsection*{Reconstructing the probability density function}
The probability density function \rev{$f_A(a)$} of a random variable \rev{$A$} is,
\rev{
\begin{subequations}
\begin{align}
        f_A(a) = \sum_{j=1}^J \Mag{ \sum_{p=0}^P A_p \frac{\dee \Phi_p}{\dee \xi}(\randomroot_j)}^{-1} W(\randomroot_j)
\intertext{where $\randomroot_j$, $j=1, \ldots, J$ are the real roots of the polynomial}
        a - \sum_{p=0}^P A_p \pcbasis_p(\xi) = 0
\end{align}\label{eqn:pdf}%
\end{subequations}}
which can be calculated numerically for a specific \rev{realisation $a$}.
Hence, the probability density function is computed by evaluating equation~\eqref{eqn:pdf} for a range of outcomes.

\subsection{Stochastic Galerkin reformulation of the deterministic model}
%In the deterministic 1D shallow water equations (equation~\ref{eqn:swe}), the flow vector is $\flow(x, t)$ and the topography is $z(x)$.
The solution of the stochastic 1D shallow water equations is now random because it depends on uncertain initial conditions, uncertain boundary conditions and uncertain topography.
Hence, the stochastic 1D shallow water equations depend not only upon space $x$ and time $t$, but additionally upon the random event $\theta$.
The stochastic flow vector $\flow(x, t, \theta)$ becomes a general stochastic process having arbitrary probability distributions that vary in space and time.
Similarly, the stochastic topography $z(x, \theta)$ has arbitrary probability distributions that vary in space.

The stochastic Galerkin reformulation of the deterministic model involves three steps to (i) replace the deterministic variables, $\flow_i^{(n)}$ and $z_i$, with random variables $\flow_i^{(n)}(\theta)$ and $z_i(\theta)$, (ii) rewrite the deterministic formulation using these random variables, and (iii) make a stochastic Galerkin projection onto the Wiener-Hermite basis.

\subsubsection*{Replacing deterministic variables with random variables}

For all elements $i=1, \ldots, M$ across all time levels, every deterministic flow variable $\flow_i^{(n)} = \left[h_i^{(n)}, q_i^{(n)}\right]^\T$ and deterministic topography variable $z_i$ becomes a random variable approximated by a Wiener-Hermite Polynomial Chaos expansion:
\begin{align}
\flow_i^{(n)}(\theta) \approx \sum _{p=0}^P \flow_{i,p}^{(n)} \pcbasis_p(\xi(\theta))
    \:\text{,}\quad
z_i(\theta) \approx \sum_{p=0}^P z_{i,p} \pcbasis_p(\xi(\theta))
\label{eqn:pc-expansion}%
\end{align}
where $\flow_{i,p}^{(n)} = \left[ h_{i,p}^{(n)}, q_{i,p}^{(n)} \right]^\T$ and $z_{i,p}$ are the $p$\textsuperscript{th} order expansion coefficients over element $i$ at time level $n$.

The reconstructed topography and reconstructed limits become functions of random variables.
The reconstructed topography at interface $i+1/2$ becomes
\begin{align}
	\sum_{p=0}^P z^\star_{i+1/2,p} \pcbasis_p
	=
	\frac{1}{2}
	\left(
	\sum_{p=0}^P z_{i,p} \pcbasis_p
	+
	\sum_{p=0}^P z_{i+1,p} \pcbasis_p
	\right)
\end{align}
and so $z^\star_{i+1/2,p} = (z_{i,p} + z_{i+1,p})/2$ due to basis orthogonality.
The reconstructed limits $\flow^{K,\star}_{i+1/2,p} = \left[ h^{K,\star}_{i+1/2,p}, q^{K,\star}_{i+1/2,p} \right]^\T$ (where $K = + \text{ or } -$) are calculated in a similar fashion.

Random variables that are functions of other random variables can be calculated in the same way.
In particular, water \rev{depth} can be expressed as a function of free-surface elevation and topography such that, due to basis orthogonality,
\begin{align}
h_{i,p}^{(n)} = \eta_{i,p}^{(n)} - z_{i,p}
\label{eqn:h-eta-z}
\end{align}
which is used later for specifying initial conditions.

\subsubsection*{Rewriting the deterministic formulation using random variables}

The deterministic finite volume formulation given by equation~\eqref{eqn:swe-discrete} is rewritten in terms of the random variables in equation~\eqref{eqn:pc-expansion}.
As a result, the numerical fluxes $\riemannflux_{i+1/2}^{(n)}$ and $\riemannflux_{i-1/2}^{(n)}$ and source term vector $\source_i^{(n)}$ become functions of random variables.
The numerical flux $\riemannflux_{i+1/2}^{(n)}$ becomes
\begin{align}
	\riemannflux_{i+1/2}^{(n)} = \riemannflux \left(
	\sum_{p=0}^P \flow^{-,\star}_{i+1/2,p} \pcbasis_p, 
	\sum_{p=0}^P \flow^{+,\star}_{i+1/2,p} \pcbasis_p
	\right)
\end{align}
and similarly for $\riemannflux_{i-1/2}^{(n)}$.
The source term vector $\source_i^{(n)}$ in equation~\eqref{eqn:source} becomes
\rev{\begin{align}
	\source_i^{(n)} = \left[ 0, -\frac{g}{\Delta x}
	\left\{
	\sum_{p=0}^P \left(\frac{h^{+,\star}_{i-1/2,p} + h^{-,\star}_{i+1/2,p}}{2} \right) \pcbasis_p \right\}
\left\{ \sum_{s=0}^P \left( z^\star_{i+1/2,s} - z^\star_{i-1/2,s} \right) \pcbasis_s \right\}
	\right]^\T
	\label{eqn:random-source}
\end{align}}
Equation~\eqref{eqn:random-source} involves the product of two expressions, each delimited by braces.
Since both expressions include Wiener-Hermite expansions then different indices, $p$ and $s$, are needed for the expansion coefficients in each expression.

\subsubsection*{Stochastic Galerkin projection}

Due to the orthogonal property (equation~\ref{eqn:orthogonal}) and distributive property (equation~\ref{eqn:distributive}) of the Wiener-Hermite basis, a Galerkin projection of equation~\eqref{eqn:swe-discrete} onto the basis functions $\pcbasis_l, l = 0, \ldots, P$ produces $P+1$ decoupled equations:
\begin{align}
    \flow_{i,l}^{(n+1)} = \flow_{i,l}^{(n)}
    - \frac{\Delta t}{\Ensemble{\pcbasis_l^2}}
    \left(
    \frac{
    \Ensemble{\riemannflux_{i+1/2}^{(n)} \pcbasis_l}
    -
    \Ensemble{\riemannflux_{i-1/2}^{(n)} \pcbasis_l}
    }{\Delta x}
    - \Ensemble{\source_i^{(n)} \pcbasis_l}
    \right) \label{eqn:swe-pc}
\end{align}
Equation~\eqref{eqn:swe-pc} involves ensemble averages of numerical fluxes, $\Ensemble{\riemannflux_{i+1/2}^{(n)} \pcbasis_l}$ and $\Ensemble{\riemannflux_{i-1/2}^{(n)} \pcbasis_l}$, and an ensemble average of the source term vector, $\Ensemble{\source_i^{(n)} \pcbasis_l}$.
There is no straightforward method for calculating an ensemble average of the numerical flux because it is nonlinear.
Instead, the integral in equation~\eqref{eqn:ensemble-average} is approximated by Gauss-Hermite quadrature,
\begin{align}
    \Ensemble{\riemannflux_{i+1/2}^{(n)} \pcbasis_l}
    \approx
    \sum_{j=1}^{P+1} w_j
    \riemannflux\left(
	\sum_{p=0}^P \flow_{i+1/2,p}^{-,\star} \pcbasis_p(\xi_j),
	\sum_{p=0}^P \flow_{i+1/2,p}^{+,\star} \pcbasis_p(\xi_j)
	\right)
    \pcbasis_l(\xi_j) W(\xi_j) \label{eqn:pc-flux}
\end{align}
where $w_j$ are the quadrature weights and $\xi_j$ are the quadrature points.
The ensemble average $\Ensemble{\riemannflux_{i-1/2}^{(n)} \pcbasis_l}$ is calculated in the same way.

Unlike the nonlinear numerical flux, the ensemble average of the source term vector $\Ensemble{\source_i^{(n)} \pcbasis_l}$ is linear and can be derived directly from equation~\eqref{eqn:random-source}:
\begin{align}
\Ensemble{\source_i^{(n)} \pcbasis_l} &= \left[ 0,
    - \frac{g}{\Delta x}
    \sum_{p=0}^P \sum_{s=0}^P
\left(\frac{h^{+,\star}_{i-1/2,p} + h^{-,\star}_{i+1/2,p}}{2}\right)
\left( z^\star_{i+1/2,s} - z^\star_{i-1/2,s} \right)
    \Ensemble{\pcbasis_p \pcbasis_s \pcbasis_l}
    \right]^\T
\label{eqn:pc-source}
\end{align}
The ensemble averages $\Ensemble{\pcbasis_p \pcbasis_s \pcbasis_l}$ in equation~\eqref{eqn:pc-source} and $\Ensemble{\pcbasis_l^2}$ in equation~\eqref{eqn:swe-pc} can be calculated analytically or exactly by Gauss-Hermite quadrature.
Since these calculations do not depend on the solution then they can be precomputed once and stored.

%% file: cproperty.tex
\subsection{Well-balanced Property with Uncertain Topography}

In this section it is shown that the stochastic Galerkin model satisfies the well-balanced property for a lake-at-rest with uncertain topography.
\rev{A lake at rest is an idealised scenario where the free-surface is entirely flat and motionless, much like still water in a lake.}
Assuming a lake-at-rest \rev{for the stochastic model}, the mean free-surface elevation is constant and the mean discharge is zero, and there is no uncertainty in the free-surface elevation or discharge.
The bed elevation can have any spatial profile and stochastic profile.

Since the free-surface elevation is constant and the averaged bed elevation is continuous at interfaces then $\eta_{i+1/2,p}^- - z^\star_{i+1/2,p} = \eta_{i+1/2,p}^+ - z^\star_{i+1/2,p}$, hence $h^{-,\star}_{i+1/2,p} = h^{+,\star}_{i+1/2,p}$ for all $i = 0, \ldots, M$ and $p = 0, \ldots, P$.
Shorthand notation is introduced such that
$\eta_p$ is the $p$\textsuperscript{th} expansion coefficient of the spatially uniform free-surface elevation, and
$h^\star_{i+1/2,p} = h^{-,\star}_{i+1/2,p} = h^{+,\star}_{i+1/2,p}$.
Since $h^\star_{i+1/2,p}$ is continuous at interfaces and the discharge is zero then the numerical flux is equal to the physical flux.

In order to preserve a lake-at-rest solution in element $i$ then the ensemble average of the flux gradient must balance the ensemble average of the source term vector in equation~\eqref{eqn:swe-pc}.
The balance of mass and momentum components can be considered separately.
The mass continuity equation balances because the discharge is zero.
For the momentum equation, the following equality must hold for all $l = 0, \ldots, P$:
\begin{align}
\frac{\Ensemble{F_{i+1/2}^{(n)} \pcbasis_l} - \Ensemble{F_{i-1/2}^{(n)} \pcbasis_l}}{\Delta x}
-
\Ensemble{S_i^{(n)} \pcbasis_l}
= 0
\label{eqn:momentum-balance-separate}
\end{align}
where the ensemble average of the bed slope source term $\Ensemble{S_i^{(n)} \pcbasis_l}$ is given by equation~\eqref{eqn:pc-source} and the momentum flux $F_{i+1/2}^{(n)}$ is
\begin{align}
F_{i+1/2}^{(n)} = \frac{g}{2}
\left\{\sum_{p=0}^P h^\star_{i+1/2,p} \pcbasis_p\right\}
\left\{\sum_{s=0}^P h^\star_{i+1/2,s} \pcbasis_s\right\}
\end{align}
and similarly for $F_{i-1/2}^{(n)}$.
Due to the distributive property of the basis (equation~\ref{eqn:distributive}), equation~\eqref{eqn:momentum-balance-separate} can be rewritten as
\begin{align}
\Ensemble{\frac{F_{i+1/2}^{(n)} 
-
F_{i-1/2}^{(n)}}{\Delta x} \pcbasis_l}
=
\Ensemble{ S_i^{(n)} \pcbasis_l}
\label{eqn:momentum-balance}
\end{align}
The ensemble average of the momentum flux gradient is
\begin{align}
    \Ensemble{\frac{F_{i+1/2}^{(n)} 
    -
    F_{i-1/2}^{(n)}}{\Delta x} \pcbasis_l}
    &=
    \frac{g}{2 \Delta x}
    \sum_{p=0}^P \sum_{\palt=0}^P
    \left( h^\star_{i+1/2,p} \, h^\star_{i+1/2,\palt} - 
    h^\star_{i-1/2,p} \, h^\star_{i-1/2,\palt} \right)
    \Ensemble{\pcbasis_p \pcbasis_\palt \pcbasis_l}
    \label{eqn:momentum-flux-unfactorised}
\intertext{To be able to factorise equation~\eqref{eqn:momentum-flux-unfactorised}, the first step is to add then subtract the term $h^\star_{i-1/2,p} h^\star_{i+1/2,s}$ to yield}
    \Ensemble{\frac{F_{i+1/2}^{(n)} 
    -
    F_{i-1/2}^{(n)}}{\Delta x} \pcbasis_l}
    &=
    \frac{g}{2 \Delta x}
    \sum_{p=0}^P \sum_{\palt=0}^P
    \left( h^\star_{i+1/2,p} \, h^\star_{i+1/2,\palt} - 
    h^\star_{i-1/2,p} \, h^\star_{i-1/2,\palt} + \right. \nonumber \\
    &\left.
    \hspace{6em}
    h^\star_{i-1/2,p} \, h^\star_{i+1/2,s} -
    h^\star_{i-1/2,p} \, h^\star_{i+1/2,s} \right)
    \Ensemble{\pcbasis_p \pcbasis_\palt \pcbasis_l}
    \label{eqn:momentum-flux-4term}
\intertext{Using the distributive property of multiplication, and associativity of summation, equation~\eqref{eqn:momentum-flux-4term} can be rewritten as}
    \Ensemble{\frac{F_{i+1/2}^{(n)} 
    -
    F_{i-1/2}^{(n)}}{\Delta x} \pcbasis_l}
    &=
    \frac{g}{2 \Delta x}
    \left[
    \sum_{p=0}^P \sum_{\palt=0}^P
    \left( h^\star_{i+1/2,p} \, h^\star_{i+1/2,\palt} - 
    h^\star_{i-1/2,p} \, h^\star_{i-1/2,\palt} + 
    h^\star_{i-1/2,p} \, h^\star_{i+1/2,s}
    \right)
    \Ensemble{\pcbasis_p \pcbasis_\palt \pcbasis_l}
    \right.
    \nonumber \\
    &\left.
    \hspace{2.5em}
    - 
    \sum_{p=0}^P \sum_{\palt=0}^P
    h^\star_{i-1/2,p} \, h^\star_{i+1/2,s}
    \Ensemble{\pcbasis_p \pcbasis_\palt \pcbasis_l} \right]
    \label{eqn:momentum-flux-split}
\intertext{
The final term in equation~\eqref{eqn:momentum-flux-split} can be rewritten as $\sum_{p=0}^P \sum_{s=0}^P h^\star_{i-1/2,s} \, h^\star_{i+1/2,p} \Ensemble{\pcbasis_p \pcbasis_s \pcbasis_l}$ due to commutativity of summation and product operators as well as the basis function (equation~\ref{eqn:commutative}).
Then, by associativity of summation:
}
    \Ensemble{\frac{F_{i+1/2}^{(n)} 
    -
    F_{i-1/2}^{(n)}}{\Delta x} \pcbasis_l}
    &=
    \frac{g}{2 \Delta x}
    \sum_{p=0}^P \sum_{\palt=0}^P
    \left( h^\star_{i+1/2,p} \, h^\star_{i+1/2,\palt} - 
    h^\star_{i-1/2,p} \, h^\star_{i-1/2,\palt} + \right. \nonumber \\
    &\left.
    \hspace{6em}
    h^\star_{i-1/2,p} \, h^\star_{i+1/2,s} -
    h^\star_{i-1/2,s} \, h^\star_{i+1/2,p} \right)
    \Ensemble{\pcbasis_p \pcbasis_\palt \pcbasis_l}
    \label{eqn:momentum-flux-4term-reordered}
\intertext{Finally, equation~\eqref{eqn:momentum-flux-4term-reordered} can be factorised as}
    \Ensemble{\frac{F_{i+1/2}^{(n)} 
    -
    F_{i-1/2}^{(n)}}{\Delta x} \pcbasis_l}
    &= 
    \frac{g}{2 \Delta x}
    \sum_{p=0}^P \sum_{\palt=0}^P
    \left( h^\star_{i-1/2,p} + h^\star_{i+1/2,p} \right)
    \left( h^\star_{i+1/2,s} - h^\star_{i-1/2,s} \right)
    \Ensemble{\pcbasis_p \pcbasis_\palt \pcbasis_l} \\
\intertext{and, given that $h^\star_{i+1/2,s} = \eta_s - \zmodified_{i+1/2,s}$ and $h^\star_{i-1/2,s} = \eta_s - \zmodified_{i-1/2,s}$, then}
    \Ensemble{\frac{F_{i+1/2}^{(n)} 
    -
    F_{i-1/2}^{(n)}}{\Delta x} \pcbasis_l}
    &=
    -\frac{g}{\Delta x}
    \sum_{p=0}^P \sum_{\palt=0}^P
    \left(
    \frac{h^{+,\star}_{i-1/2,p} + h^{-,\star}_{i+1/2,p}}{2}
    \right)
    \left( \zmodified_{i+1/2,\palt} - \zmodified_{i-1/2,\palt} \right)
    \Ensemble{\pcbasis_p \pcbasis_\palt \pcbasis_l}
    \label{eqn:momentum-flux}
\end{align}
Equation~\eqref{eqn:momentum-flux} is equal to the ensemble average of the bed slope source term $\Ensemble{S_i^{(n)} \pcbasis_l}$ given in equation~\eqref{eqn:pc-source}, hence discrete balance is preserved.

%% file: experiments.tex
\section{\rev{Idealised numerical tests}}

\rev{Three} numerical \rev{tests} are performed to validate the well-balanced stochastic Galerkin model.
The first test simulates a lake-at-rest over \rev{idealised,} uncertain topography to verify that the stochastic Galerkin model preserves well-balancing numerically.
The second test simulates flow over a hump with an uncertain elevation.
This test is designed to challenge the stochastic Galerkin model at representing discontinuous, non-Gaussian probability distributions by generating a steady-state solution that may be subcritical or transcritical depending on the hump elevation.
\rev{The third test verifies the robustness of the stochastic Galerkin model for a steady flow over a highly irregular and uncertain bed that is more representative of real-world river hydraulics.}

\subsection{\rev{Specification of Tests 1 and 2}}
\rev{For the first two tests,} the 1D domain is [\SI{-50}{\meter}, \SI{50}{\meter}], tessellated by $M = 100$ elements with no overlaps or gaps such that the mesh spacing is $\Delta x = \SI{1}{\meter}$.
The timestep is $\Delta t = \SI{0.15}{\second}$ resulting in a maximum Courant number of about $0.8$.
By choosing a fixed timestep, simulations of a given test complete in the same number of timesteps irrespective of the model configuration, and error accumulation due to timestepping errors will be the same across all models.

Both tests include a topographic hump centred at $x = \SI{0}{\meter}$ with a region of Gaussian uncertainty.
Following a similar approach to \citet{ge2008}, there are two representations for the same uncertain topography.
The first representation enables smooth topography profiles to be randomly generated in Monte Carlo iterations.
The topography $z$ is defined as
\rev{\begin{align}
z(x, \hump) = \hump \sech^2 \left( \frac{\pi x}{\lambda} \right) \label{eqn:bed}
\end{align}}
where the hump amplitude \rev{$\hump$} is a random variable with mean \rev{$\humpmean = \SI{0.6}{\meter}$} and standard deviation \rev{$\sigma_\hump = \SI{0.3}{\meter}$}, and the half-width is $\lambda = \SI{10}{\meter}$.
\rev{This topography profile is seen in figure~\ref{fig:criticalSteadyState-flow}.}
The second representation is used by the stochastic Galerkin model, with topography represented by expansion coefficients $z_{i,0}, \ldots, z_{i,P}$.
To be able to calculate the topography expansion coefficients, equation~\eqref{eqn:bed} must be expressed in terms of mean topography $\zmean(x)$ and topographic variance $\sigma_z^2(x)$ without involving the random variable \rev{$\hump$}.
The mean topography $\zmean(x)$ is simply
\rev{
\begin{align}
    \zmean(x) = z(x, \humpmean)
    \label{eqn:bed-mean}
\end{align}}
The topographic variance is
\rev{\begin{align}
    \sigma^2_z(x) = \E[z^2(x, \hump)] - \E^2[z(x, \hump)]
    \label{eqn:bed-variance}
\end{align}}
Equation~\eqref{eqn:bed-variance} can be rewritten using Taylor series expansions of the two terms \rev{$\E\left[z^2(x, \hump)\right]$} and \rev{$\E^2\left[z(x, \hump)\right]$}.
To illustrate the approach, a Taylor series expansion of \rev{$\E\left[ z(x, \hump) \right]$} about \rev{$\humpmean$} is 
\begin{align}
    \E\left[ z(x, \hump) \right] &= \E\left[ z(x, \humpmean + (\hump - \humpmean)) \right] \nonumber \\
    &= \E\left[ z + \frac{\partial z}{\partial \hump} (\hump - \humpmean) + \frac{1}{2} \frac{\partial^2 z}{\partial \hump^2} \left(\hump - \humpmean\right)^2 + \mathcal{O}(\hump^3) \right]
    \label{eqn:mean-taylor}
\intertext{where $z$ is shorthand for \rev{$z(x, \humpmean)$} and \rev{$\mathcal{O}(\hump^3)$} is the error term involving high-order derivatives \rev{$\partial^m z/\partial \hump^m$} with $m \geq 3$.
Since \rev{$\E\left[ \hump - \humpmean \right] = 0$} and \rev{$\E\left[ \left(\hump-\humpmean\right)^2\right] = \sigma_\hump^2$} then equation~\eqref{eqn:mean-taylor} simplifies to}
    \E\left[ z(x, \hump) \right] &= z + \frac{1}{2}\frac{\partial^2 z}{\partial \hump^2} \sigma_\hump^2 + \mathcal{O}(\hump^3)
\end{align}
Applying this approach to equation~\eqref{eqn:bed-variance} gives:
\begin{align}
    \sigma_z^2(x) &=
    z^2 +
    \left[
    \left(\frac{\partial z}{\partial \hump}\right)^2
    + z \frac{\partial^2 z}{\partial \hump^2}
    \right]
    \sigma^2_\hump
    -
    \left[
    z + \frac{1}{2} \frac{\partial^2 z}{\partial \hump^2} \sigma_\hump^2
    \right]^2 + \mathcal{O}(\hump^3) \label{eqn:z-taylor}
\intertext{For the topographic profile given by equation~\eqref{eqn:bed}, it holds that \rev{$\partial^m z/\partial \hump^m = 0$} where $m \geq 2$, \rev{so the Taylor series approximation introduces no spurious oscillations in stochastic space.
The topographic variance in equation~\eqref{eqn:z-taylor} then simplifies to}}
    \sigma_z^2(x) &= \left( \frac{\partial z(x, \humpmean)}{\partial \hump} \sigma_\hump \right)^2 \label{eqn:z-variance-continuous}
\end{align}
Equipped with analytic expressions for the mean topography $\zmean(x)$ (equation~\ref{eqn:bed-mean}) and topographic variance $\sigma_z^2(x)$ (equation~\ref{eqn:z-variance-continuous}), now the topography expansion coefficients $z_{i,0}, \ldots, z_{i,P}$ can be calculated.
Since the topographic bump has a Gaussian probability distribution with $\mu_1[z(x)] = \zmean(x)$, $\mu_2[z(x)] = \sigma_z^2(x)$ and high-order moments $\mu_m[z(x)] = 0$ for $m \geq 3$ then, using equations~\eqref{eqn:moment}, \eqref{eqn:mean} and \eqref{eqn:variance}, the topography expansion coefficients are
\begin{align}
    z_{i,p} = \begin{cases}
    \zmean(x_i) & \text{if $p=0$} \\
    \sigma_z(x_i) & \text{if $p=1$} \\
    0 & \text{otherwise}
    \end{cases}
    \label{eqn:z-pc-coeffs}
\end{align}
where values are calculated at the centre point $x_i$ for all elements $i=1,\ldots, M$.

%and $z_{i,1}$ is calculated using equations~\eqref{eqn:z-variance-continuous} and \eqref{eqn:variance} with $\mu_2[z(x, \xi)] = \sigma_z^2(x)$.

The initial water \rev{depth} expansion coefficients $h_{i,0}^{(0)}, \ldots, h_{i,P}^{(0)}$ can be calculated in terms of free-surface elevation and topography using equation~\eqref{eqn:h-eta-z}.
For both tests, the initial, spatially-uniform mean free-surface elevation is \SI{1.5}{\meter} with no initial uncertainty such that $\eta_{i,0}^{(0)} = \SI{1.5}{\meter}$ and $\eta_{i,p}^{(0)} = 0$ with $p = 1, \ldots, P$ and $i = 1, \ldots, M$.

\input{lakeAtRest}

\input{criticalSteadyState}
\input{tsengSteadyState}

%% file: lakeAtRest.tex
\subsection{\rev{Test 1:} Lake-at-Rest Over an Uncertain Bed}

\rev{Numerical methods that are not well-balanced produce spurious waves in the vicinity of sloping topography, and these spurious waves are particularly evident for slow-moving flows with weak momentum fluxes.
In the limit, the momentum flux is zero and the water is motionless, like the water surface on a lake at rest.
Hence, the lake-at-rest test is ideally suited to verify the well-balanced property, since the analytic solution preserves the resting state forever, and any waves generated by a numerical model are entirely spurious \citep{bermudez-vazquez1994}.}

To present a challenging test, a rectangular obstacle is introduced to the right of the uncertain hump given by equation~\eqref{eqn:bed}, so the bed elevation $z$ becomes
\begin{subequations}
\begin{align}
    z(x, \hump) &= \hump \sech^2 \left( \frac{\pi x}{\lambda} \right) + \zrectbump(x) \text{,}
    \intertext{where $\zrectbump$ is the rectangular obstacle:}
    \zrectbump(x) &= \begin{cases}
    \humpmean & \text{if $30 < x \leq 40$,} \\
    0 & \text{otherwise}
    \end{cases}
\end{align}
\end{subequations}

Results of the well-balanced stochastic model are compared with those of a stochastic model having a centred difference approximation of the source term vector that does not exactly balance the numerical flux gradient.
The centred difference model is the same as the well-balanced stochastic model except for two changes.
First, numerical fluxes are calculated using the original, unmodified flow variables:
\begin{align}
	\riemannflux_{i+1/2}^{(n)} = \riemannflux \left(
	\sum_{p=0}^P \flow^-_{i+1/2,p} \pcbasis_p, 
	\sum_{p=0}^P \flow^+_{i+1/2,p} \pcbasis_p
	\right) \label{eqn:flux-centred}
\end{align}
Second, the ensemble average of the source term vector uses a centred difference approximation:
\rev{\begin{align}
    \Ensemble{\source_i \pcbasis_l} =
    \left[ 0, -g \sum_{p=0}^P \sum_{s=0}^P h_{i,p}
    \frac{z_{i+1,s} - z_{i-1,s}}{2 \Delta x}
    \Ensemble{\pcbasis_p \pcbasis_s \pcbasis_l} \right]^\T \label{eqn:source-centred}
\end{align}}
The centred difference model and the well-balanced model are both configured with a Wiener-Hermite basis order $P = 3$.

\begin{figure}
\centering
\begin{subfigure}{\textwidth}
\phantomsubcaption\label{fig:lakeatrest:centred:eta}
\phantomsubcaption\label{fig:lakeatrest:sgm:eta}
\phantomsubcaption\label{fig:lakeatrest:centred:q}
\phantomsubcaption\label{fig:lakeatrest:sgm:q}
\centering
\includegraphics{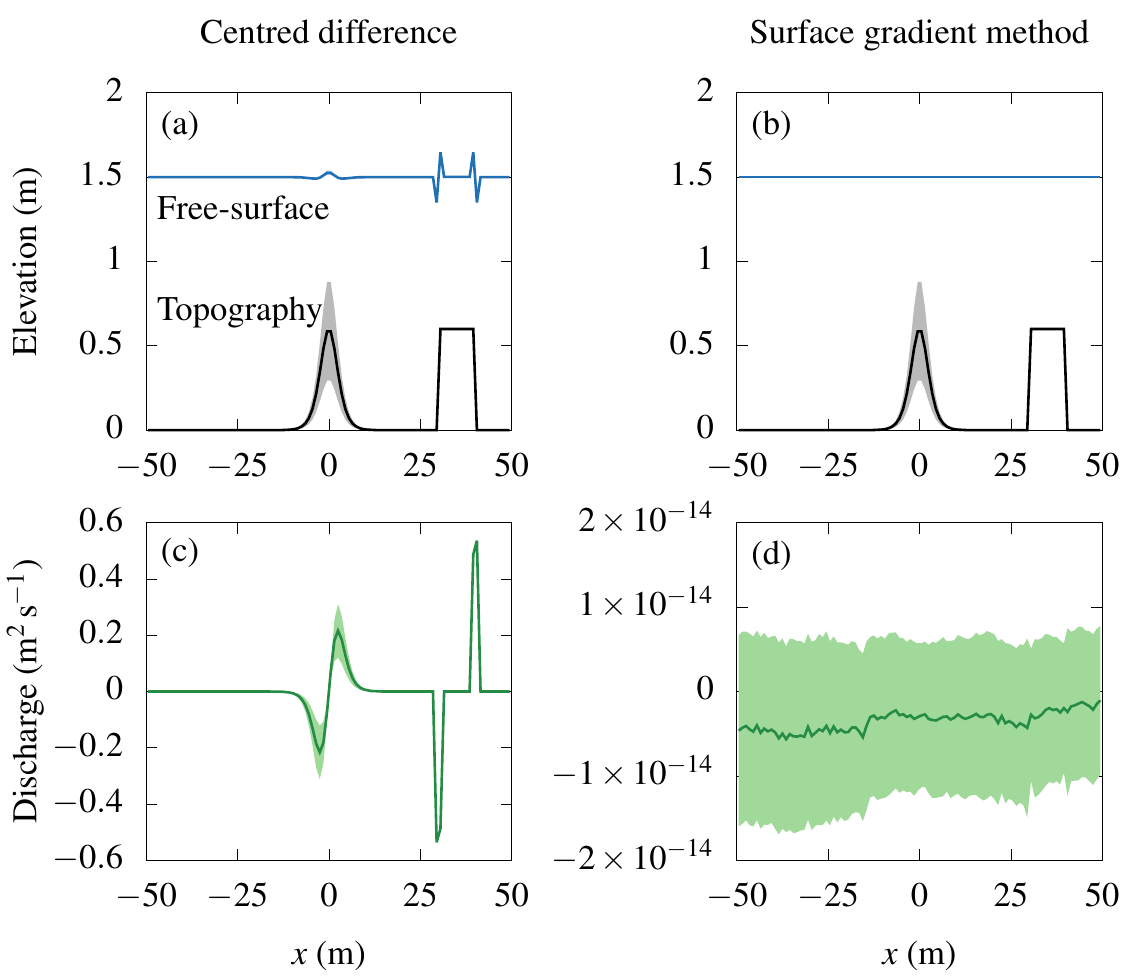}
\end{subfigure}
\caption{Stochastic lake-at-rest solutions at $t = \SI{100}{\second}$.
Mean values are marked by solid lines and shaded regions represent one standard deviation from the mean.}
\label{fig:lakeatrest}
\end{figure}

The simulated time is \SI{100}{\second} corresponding to about 670 timesteps, and the solutions for the centred difference and well-balanced stochastic Galerkin models are shown in figure~\ref{fig:lakeatrest}.
The lack of well-balancing is apparent using the centred difference model: grid-scale standing waves develop at the discontinuities either side of the rectangular hump (figure~\ref{fig:lakeatrest:centred:eta}, \ref{fig:lakeatrest:centred:q}), and a smooth standing wave also develops over the uncertain hump.
These errors persist throughout the simulation.
In contrast, the well-balanced stochastic Galerkin model preserves the initial resting state with discharges accurate to machine precision (figure~\ref{fig:lakeatrest:sgm:q}).
This numerical result confirms that the stochastic Galerkin model is well-balanced in theory and in practice.

The choice of the Wiener-Hermite basis introduces a particular limitation that imposes an upper bound on the basis order $P$, and constrains the minimum water depth that the stochastic Galerkin model can represent.
This limitation arises because the hump amplitude has a Gaussian probability distribution so the tails of the distribution extend to $\pm \infty$, meaning that there is a non-zero probability that the water depth is negative.
The stochastic Galerkin formulation presented here does not accommodate wetting-and-drying processes, and any negative water depth will crash the model.
If the basis order $P$ is increased then the Gauss-Hermite quadrature points in equation~\eqref{eqn:pc-flux} extend further into the tails of the probability distributions, leading to negative water depths being provided as input to the Riemann solver.
Similarly, raising the topography, decreasing the initial water depth, or increasing the topographic uncertainty can all produce negative water depths in the stochastic Galerkin model.
This behaviour has been verified experimentally by varying the model basis order, initial conditions and topography profile.

%% file: criticalSteadyState.tex
\subsection{\rev{Test 2:} Steady-State Critical Flow Over an Uncertain Bed}

This test is designed to challenge the stochastic Galerkin model at representing highly non-Gaussian and discontinuous distributions of stochastic flow.
The uncertain topography profile and inflow boundary condition are chosen in order to produce a nonlinear response such that the steady-state solution may be subcritical or transcritical depending on the bed elevation.
Results of the stochastic Galerkin model are validated against a Monte Carlo simulation that serves as a reference solution.

\subsubsection{Setup to Produce a Nonlinear Flow Response}
The uncertain topography is the smooth hump given by equation~\eqref{eqn:z-pc-coeffs}.
Subcritical boundary conditions are imposed such that the mean upstream discharge \rev{per unit-width} is \SI{1.65}{\meter\squared\per\second} and the mean downstream water \rev{depth} is \SI{1.5}{\meter}, with no uncertainty for the upstream discharge and downstream water \rev{depth}.
Transmissive boundary conditions are used for the upstream water \rev{depth} and downstream discharge.
These boundary conditions are chosen so that the flow is exactly critical for the mean hump amplitude $\humpmean = \SI{0.6}{\meter}$ at $x = \SI{0}{\meter}$.
Since the hump amplitude is uncertain then the flow regime is also uncertain: if the hump amplitude is less than $\humpmean$ then the flow remains subcritical; if the hump amplitude is greater than $\humpmean$ then the flow regime becomes transcritical.
\rev{In the transcritical regime, the flow upstream of the hump is subcritical, transitioning to supercritical flow over the hump.
A hydraulic jump occurs on the downstream side of the hump, where the flow becomes subcritical once more.}

\begin{figure}
    \centering
    \includegraphics{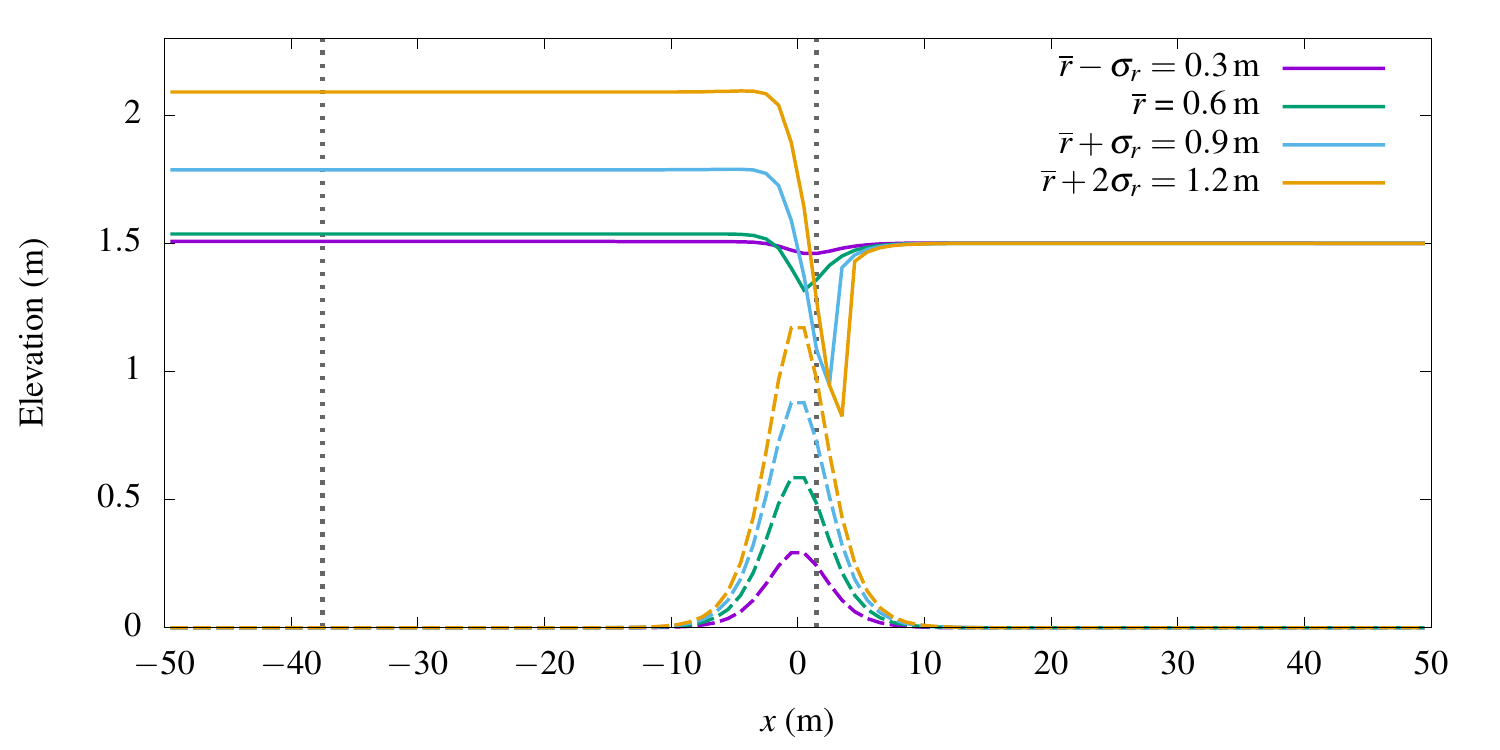}
    \caption{Well-balanced deterministic solutions of steady-state flow at $t = \SI{500}{\second}$ using four hump amplitudes, \rev{$\humpmean - \sigma_\hump = \SI{0.3}{\meter}, \humpmean = \SI{0.6}{\meter}, \humpmean+\sigma_a = \SI{0.9}{\meter}$ and $\humpmean + 2\sigma_\hump = \SI{1.2}{\meter}$}.
    The free-surface elevation is shown with solid lines and the topography profile is shown with dashed lines.
    Vertical dotted lines at $x = \SI{-37.5}{\meter}$ and $x = \SI{1.5}{\meter}$ mark the positions of the probability densities shown in figure~\ref{fig:criticalSteadyState-pdf}.}
    \label{fig:criticalSteadyState-examples}
\end{figure}

To illustrate this change in flow regime, figure~\ref{fig:criticalSteadyState-examples} shows four deterministic solutions using four different hump amplitudes.
Solutions from the well-balanced deterministic model are obtained at $t = \SI{500}{\second}$ when the water has converged on a steady state.
Convergence is measured by calculating the $L^2$ difference in mean water \rev{depth} between the current and previous timesteps:
\begin{align}
    L^2 \text{ difference in mean water \rev{depth}} = \sqrt{\sum_{i=1}^M \left(h_{i,0}^{(n)} - h_{i,0}^{(n-1)}\right)^2} \label{eqn:convergence}
\end{align}
By $t = \SI{500}{\second}$ all four deterministic solutions have converged down to \rev{a convergence error of} $10^{-4}$ \si{\meter}.
For a small hump with amplitude \rev{$\humpmean - \sigma_\hump = \SI{0.3}{\meter}$}, the flow remains subcritical.
A linear increase in hump amplitude produces a strongly nonlinear response in the steady-state water profile, as seen in figure~\ref{fig:criticalSteadyState-examples}. 
Two nonlinear responses are evident in particular: first, the upstream boundary condition allows the upstream water \rev{depth} to increase nonlinearly; second, a transcritical shock develops over the hump that increases in amplitude and moves further downstream with larger hump amplitudes.
Downstream of the hump, the water \rev{depth} is \SI{1.5}{\meter} irrespective of the hump amplitude, with this profile having propagated upstream from the imposed downstream boundary.

\subsubsection{Configuration of the Monte Carlo Reference Simulation}
The Monte Carlo reference simulation is performed using iterations of the well-balanced deterministic model.
\rev{It is necessary to perform a sufficient number of Monte Carlo iterations to ensure that the flow statistics are accurate, and more iterations are needed for more complex probability distributions.
After each Monte Carlo iteration, the mean and standard deviation of water depth at $x = \SI{1.5}{\meter}$ are measured, where the probability distribution is most complex.
The statistics are compared with those of the previous iteration, and the Monte Carlo simulation terminates once the statistics change negligibly between iterations, when statistical convergence is achieved.
Guided by these measurements, two thousand Monte Carlo iterations are needed to achieve statistical convergence for this test.}

For each Monte Carlo iteration, the topography is randomly generated using a hump amplitude drawn from the Gaussian distribution given by \rev{$(\humpmean, \sigma_\hump)$} and so the topography will always be smooth.
If instead the topography was randomly generated using $(\zmean(x), \sigma_z(x))$ then randomisation would be different in every element and the topography would not be smooth, so many more iterations would be needed to sample the stochastic solution space.
For the Monte Carlo iterations, the hump amplitude \rev{$\hump$} is constrained such that \rev{$\SI{0}{\meter} \leq \hump \leq \SI{1.4}{\meter}$} to avoid negative water depths.

\subsubsection{Spatial Profiles of the Uncertain Free-Surface Elevation}

\begin{figure}
    \centering
    \includegraphics{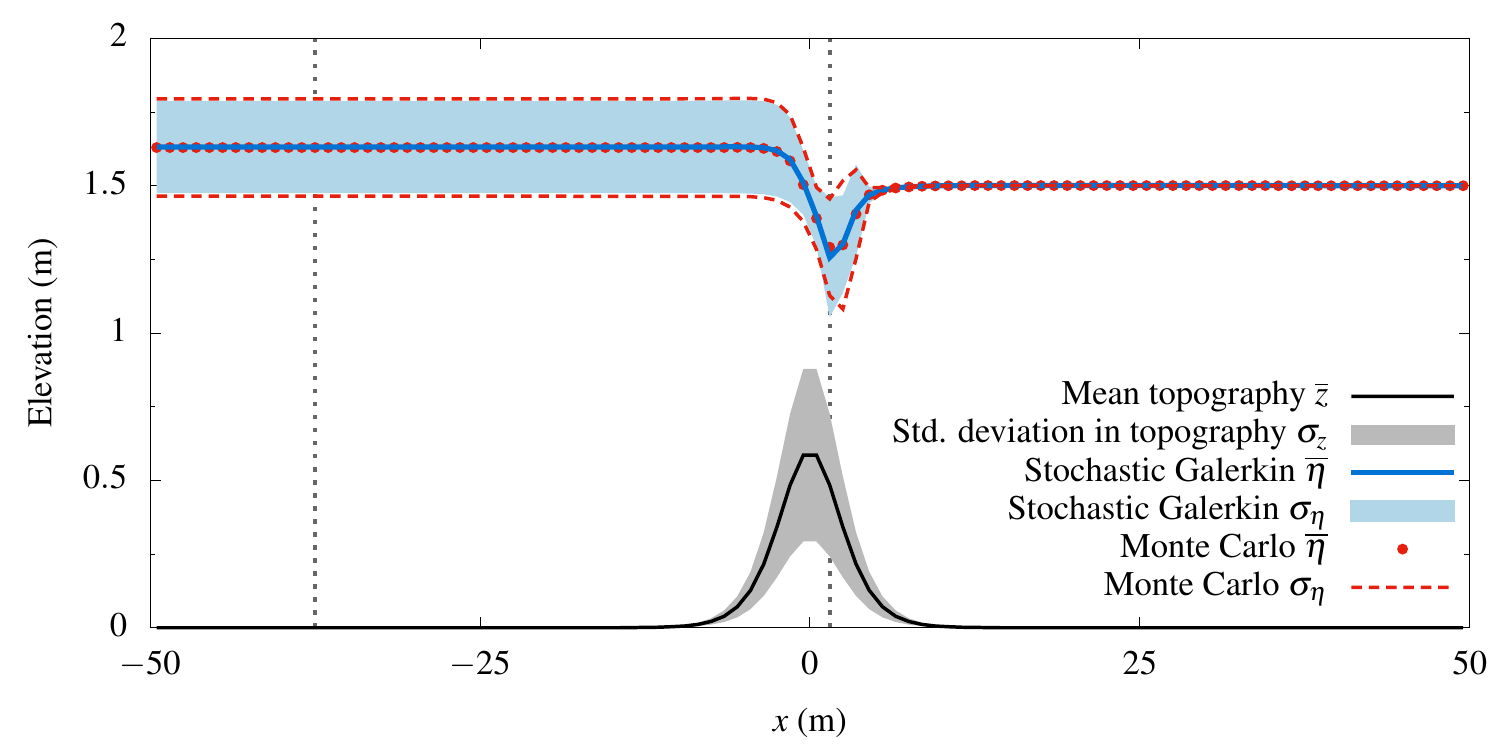}
    \caption{Solutions of steady state critical flow over an uncertain hump at $t = \SI{500}{\second}$, comparing stochastic Galerkin and Monte Carlo profiles of mean free-surface elevation $\etamean$ and standard deviation of free-surface elevation $\sigma_\eta$.
    The stochastic Galerkin result is obtained using basis order $P = 3$.
    Vertical dotted lines at $x = \SI{-37.5}{\meter}$ and $x = \SI{1.5}{\meter}$ mark the positions of the probability densities shown in figure~\ref{fig:criticalSteadyState-pdf}.
    }
    \label{fig:criticalSteadyState-flow}
\end{figure}

In figure~\ref{fig:criticalSteadyState-flow}, spatial profiles of the uncertain free-surface elevation are obtained at $t = \SI{500}{\second}$ when \rev{the water depth profiles from the Monte Carlo and stochastic Galerkin simulations have converged down to $10^{-4}$ \si{\meter} as defined by equation~\eqref{eqn:convergence}}.
Using a Wiener-Hermite basis of order $P=3$, the stochastic Galerkin model accurately represents the Monte Carlo reference profiles of the mean and standard deviation of free-surface elevation.
Upstream of the hump, the stochastic Galerkin model predicts a standard deviation that is slightly too small compared to the Monte Carlo reference solution.
Small errors in the stochastic Galerkin free-surface elevation are also visible above the hump where the flow is most complex.
Similar results are obtained using basis order $P=1$ or $P=2$, with stochastic Galerkin errors increasing slightly as the basis order is decreased (not shown).

\subsubsection{Monte Carlo Histograms of Uncertain Free-Surface Elevation}
The mean and standard deviation statistics are useful for summarising the spatial profile of uncertainty, but they are less meaningful for non-Gaussian probability distributions.
In such cases, it is more meaningful to study the complete probability distributions, which is also particularly important for flood risk assessments that are concerned with extreme events that occur in the tails of the distributions \citep{ge2011}.

\begin{figure}
    \centering
    \begin{subfigure}{\textwidth}
    \phantomsubcaption\label{fig:criticalSteadyState-pdf:upstream}
    \phantomsubcaption\label{fig:criticalSteadyState-pdf:downstream}
    \centering
    \includegraphics{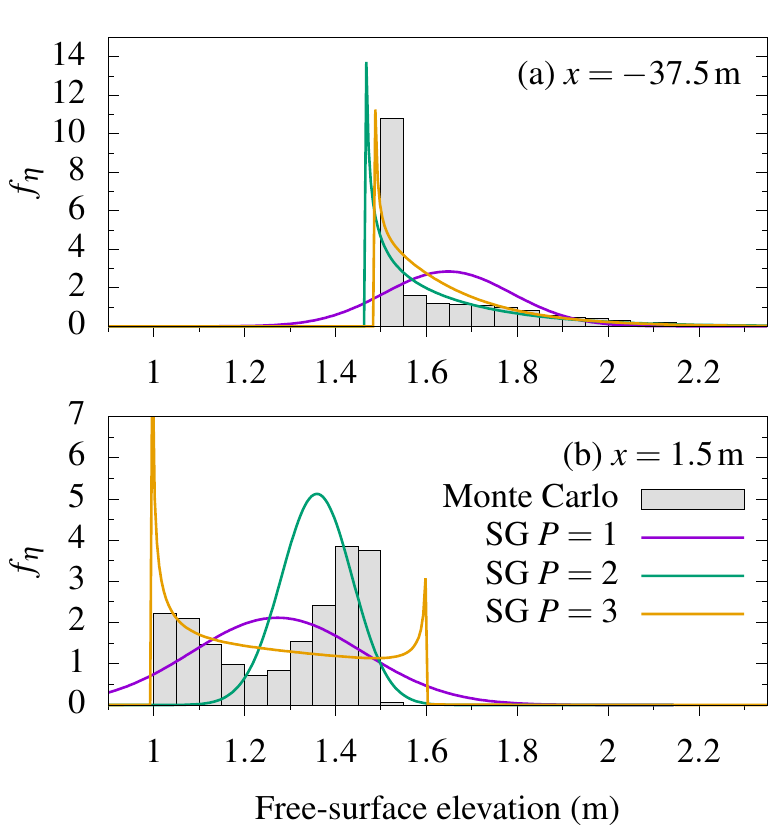}
    \end{subfigure}
    \caption{Probability distributions of free-surface elevation, $f_\eta$, at (a) $x = \SI{-37.5}{\meter}$ and (b) $x = \SI{1.5}{\meter}$ for steady state critical flow over an uncertain hump at $t = \SI{500}{\second}$.
    For the Monte Carlo reference simulation, probability distributions are estimated by histograms.
    Continuous probability density functions are reconstructed from stochastic Galerkin (SG) results using basis orders $P=1$, $P=2$ and $P=3$.}
    \label{fig:criticalSteadyState-pdf}
\end{figure}

Probability distributions of the free-surface elevation are sampled at two points at $t =
\SI{500}{\second}$: the first at $x =
\SI{-37.5}{\meter}$ and the second at $x
= \SI{1.5}{\meter}$, with these positions marked by dotted lines in
figure~\ref{fig:criticalSteadyState-examples} and
figure~\ref{fig:criticalSteadyState-flow}.
The first point is far upstream of the hump where the free-surface elevation is uncertain and spatially uniform.
The second point is immediately downstream of the hump in the region where transcritical shocks can develop.

Figure~\ref{fig:criticalSteadyState-pdf} shows histograms from the Monte Carlo reference simulation that estimate the true probability densities at the two points.
Stochastic Galerkin results which appear on the same figure are discussed later.
For each of the two points, water \rev{depths} from the 2000 Monte Carlo iterations are binned into intervals of \SI{0.05}{\meter}, and the magnitude of each bin represents the probability that the water \rev{depth} is within the given interval.
Since the histogram estimates the probability density then the total shaded area over all bins is equal to one.

The Monte Carlo histogram at $x = \SI{-37.5}{\meter}$ is shown in figure~\ref{fig:criticalSteadyState-pdf:upstream} and is discussed first.
For subcritical flows over small humps, the upstream water level remains at its initial height of \SI{1.5}{\meter}.
For transcritical flows over larger humps, the upstream water level increases nonlinearly.
Since the initial conditions and boundary conditions are chosen so that the mean flow is critical, then about 50\% of the flows are subcritical, resulting in a large peak in the
probability distribution at $\eta = \SI{1.5}{\meter}$.
The other 50\% of the flows are transcritical with elevated upstream water levels, resulting in a long tail in the distribution.

The stochastic flow response immediately downstream of the hump at $x = \SI{1.5}{\meter}$ is more complex.
For subcritical flows, the steady-state free-surface elevation at this point will be only slightly lower than the initial free-surface elevation of \SI{1.5}{\meter}.
For transcritical flows, the steady-state free-surface elevation may be much lower since,  at $x = \SI{1.5}{\meter}$, the transcritical shock is close to its minimum depth.
The subcritical and transcritical regimes appear as a bimodal distribution in the histogram in figure~\ref{fig:criticalSteadyState-pdf:downstream} with one peak around $\eta = \SI{1.0}{\meter}$ associated with transcritical flow and a second peak around $\eta = \SI{1.5}{\meter}$ associated with subcritical flow.

\subsubsection{Stochastic Galerkin Free-Surface Elevation Probability Densities}
Overlayed on the Monte Carlo histograms, figure~\ref{fig:criticalSteadyState-pdf} also shows probability density functions obtained from stochastic Galerkin simulations.
Three stochastic Galerkin simulations are performed using Wiener-Hermite bases of order $P=1$, $P=2$ and $P=3$.
Free-surface elevation expansion coefficients are calculated by rearranging equation~\eqref{eqn:h-eta-z}, from which probability density functions are reconstructed using equation~\eqref{eqn:pdf}.
If a probability density function was in exact agreement with a Monte Carlo histogram then the line would pass through the top of every histogram bin, and so any deviation represents a numerical error associated with the stochastic Galerkin model.

Using basis order $P=1$, the stochastic Galerkin model can only represent Gaussian distributions, with its two expansion coefficients representing the mean and standard deviation.
While figure~\ref{fig:criticalSteadyState-flow} confirms that the mean and standard deviation statistics from the stochastic Galerkin model are accurate, the Monte Carlo histograms cannot be well-represented by Gaussian distributions.

Increasing the basis order to $P=2$, the probability distribution of upstream water levels is in good agreement with the Monte Carlo histogram (figure~\ref{fig:criticalSteadyState-pdf:upstream}), though an error is noticeable at the discontinuity:
the Monte Carlo histogram has a discontinuity at $\eta = \SI{1.5}{\meter}$ because the upstream water level can only rise above its initial height, so the probability that $\eta < \SI{1.5}{\meter}$ is zero.
The stochastic Galerkin model with basis order $P = 2$ slightly underestimates this discontinuity, placing it at about $\eta = \SI{1.45}{\meter}$.
The distribution in free-surface elevation at $x = \SI{1.5}{\meter}$ (figure~\ref{fig:criticalSteadyState-pdf:downstream}) is not improved using basis order $P=2$, with the distribution remaining close to Gaussian.

Finally, the basis order is increased to $P=3$.
The distribution of upstream water levels shifts slightly to the right in figure~\ref{fig:criticalSteadyState-pdf:upstream}, corresponding to slightly higher water levels.
As a result, the discontinuity is closer to the true value of $\eta = \SI{1.5}{\meter}$, but this shift also produces slightly larger errors in the tail of the distribution.
The most notable improvement is seen in the distribution at $x = \SI{1.5}{\meter}$ (figure~\ref{fig:criticalSteadyState-pdf:downstream}):
the distribution from the stochastic Galerkin model now has a similar shape to the Monte Carlo histogram.
The lower bound at $\eta = \SI{1.0}{\meter}$ is accurately represented, but the true upper bound around $\eta = \SI{1.5}{\meter}$ is overestimated by about \SI{0.1}{\meter}.
The probability density function has singularities at the lower and upper bounds that appear as spikes in the plot.
These two singularities occur because the probability density function is the derivative of the cumulative density function, which is discontinuous and non-differentiable at these points.

While the stochastic Galerkin model with basis order $P=3$ adequately represented the true distribution bounds and distribution shapes, the probability density functions were not entirely accurate.
Such inaccuracies are to be expected because low-order Wiener-Hermite bases cannot represent complex distributions.
The basis order cannot be increased beyond $P=3$ because water depths are small near transcritical shocks, and the stochastic Galerkin model crashes for the same reason as discussed in the lake-at-rest test.
Instead, results might be improved by choosing a more sophisticated method to discretise stochastic space.
One candidate is the stochastic Galerkin multiwavelet approach by \citet{pettersson2014}, which is able to simulate stochastic gas dynamics with densities close to zero,  analogous to very small water depths in shallow water flows.

\subsubsection{Storage Requirements and Computation Time}

Since the stochastic Galerkin model is necessarily more complex than the deterministic model, associated increases in storage requirements and computation time are expected.
Stochastic Galerkin storage requirements scale linearly with the chosen basis order $P$ because the model stores $P+1$ expansion coefficients per variable per element.
The ensemble averages $\Ensemble{\pcbasis_p \pcbasis_s \pcbasis_l}$ in equation~\eqref{eqn:pc-source} and $\Ensemble{\pcbasis_l^2}$ in equation~\eqref{eqn:swe-pc} are constant values that can be precomputed once and stored.
Due to the commutative property of the ensemble average of basis function products (equation~\ref{eqn:commutative}), the associated storage requirements are small.

The expected increase in computation time can be estimated by examining the model formulation and assuming that calculations are performed without parallelisation.
Choosing basis order $P$, the entire evolution equation (equation~\ref{eqn:swe-pc}) must be evaluated $P+1$ times for the basis functions $0, \ldots, P$.
For each evaluation, the ensemble average of numerical fluxes is calculated by sampling the Riemann solver $P+1$ times using Gauss-Hermite quadrature (equation~\ref{eqn:pc-flux}).
The computation time for the ensemble average of the source term vector (equation~\ref{eqn:pc-source}) can be neglected because its calculation is trivially fast compared to the Riemann solver. 
Hence, by considering the total samples of the Riemann solver, it is expected that the stochastic Galerkin model will be at least $\left(P+1\right)^2$ times slower than the deterministic solver.

Measuring the elapsed CPU time to simulate the steady-state test confirms this expectation: the stochastic Galerkin model with basis order $P=3$ is in fact about 20 times slower than the deterministic model, and the Riemann solver accounts for about 90\% of the stochastic Galerkin total computation time.
Compared to the Monte Carlo simulation which uses \num{2000} deterministic iterations, the stochastic Galerkin model is about 100 times faster.
Note that the deterministic model, stochastic Galerkin model and Monte Carlo simulation were implemented without parallelisation.

The speed-up observed here compares favourably with numerical \rev{tests} performed by \citet{ge2008}.
Their stochastic Galerkin model with basis order $P=5$ was about 200 times slower than their deterministic model.
Their Monte Carlo simulation used \num{10000} iterations, making their stochastic Galerkin model about 50 times faster.
A direct comparison of computation time is not attempted since \citet{ge2008} used a second-order finite volume formulation that is inherently more expensive than the first-order formulation presented here.

\subsubsection{Opportunities for Parallel Computation}

Computation time for Monte Carlo and stochastic Galerkin simulations could be reduced by exploiting opportunities for parallelism.
Monte Carlo simulations are called `embarrassingly parallel' because, given sufficient processors, it is easy to perform iterations entirely in parallel.
The stochastic Galerkin model is not embarrassingly parallel, but most operations can be parallelised nevertheless.
Due to basis orthogonality, the evolution equation (equation~\ref{eqn:swe-pc}) can be evaluated in parallel over basis functions $\pcbasis_0, \ldots, \pcbasis_P$.
Gauss-Hermite quadrature of numerical fluxes (equation~\ref{eqn:pc-flux}) can be parallelised, too.
In this way, sampling the Riemann solver, which accounts for about 90\% of the total computational cost, could be fully parallelised.

In theory, by fully exploiting parallelism, Monte Carlo and stochastic Galerkin simulations could be made to run almost as fast as a single iteration of the deterministic model, but this assumes access to sufficient processors.
In practice, a fully parallelised Monte Carlo would require thousands of processors which are unavailable on typical hardware.
In contrast, the stochastic Galerkin model with basis order $P=3$ could be fully parallelised with just 16 processors using current, commodity hardware.
With more hardware, additional parallelism could be achieved using domain decomposition techniques.
The stochastic Galerkin method imposes no barriers to domain decomposition because it preserves the local element-wise operations of the underlying deterministic formulation.
Given these substantial reductions in computation time and hardware demands, a stochastic Galerkin shallow flow model could alleviate the computational constraints associated with conventional Monte Carlo simulations \citep{neal2013}, and allow probabilistic simulations to account for more sources of uncertainty.

%% file: tsengSteadyState.tex
\rev{\subsection{Test 3: Steady-State Slow Flow Over an Irregular Bed}}
\rev{While the previous tests are restricted to idealised topography profiles, this final test simulates slow flow over a highly irregular bed that is more representative of real-world river hydraulics.
The bed is defined with a range of uncertainty that characterises bathymetric measurement errors.
Stochastic Galerkin model results are verified against the analytic solution.}

\begin{figure}
    \centering
    \begin{subfigure}{\textwidth}
    \centering
    \phantomsubcaption\label{fig:tsengSteadyState-flow:elevation}
    \phantomsubcaption\label{fig:tsengSteadyState-flow:u}
    \includegraphics{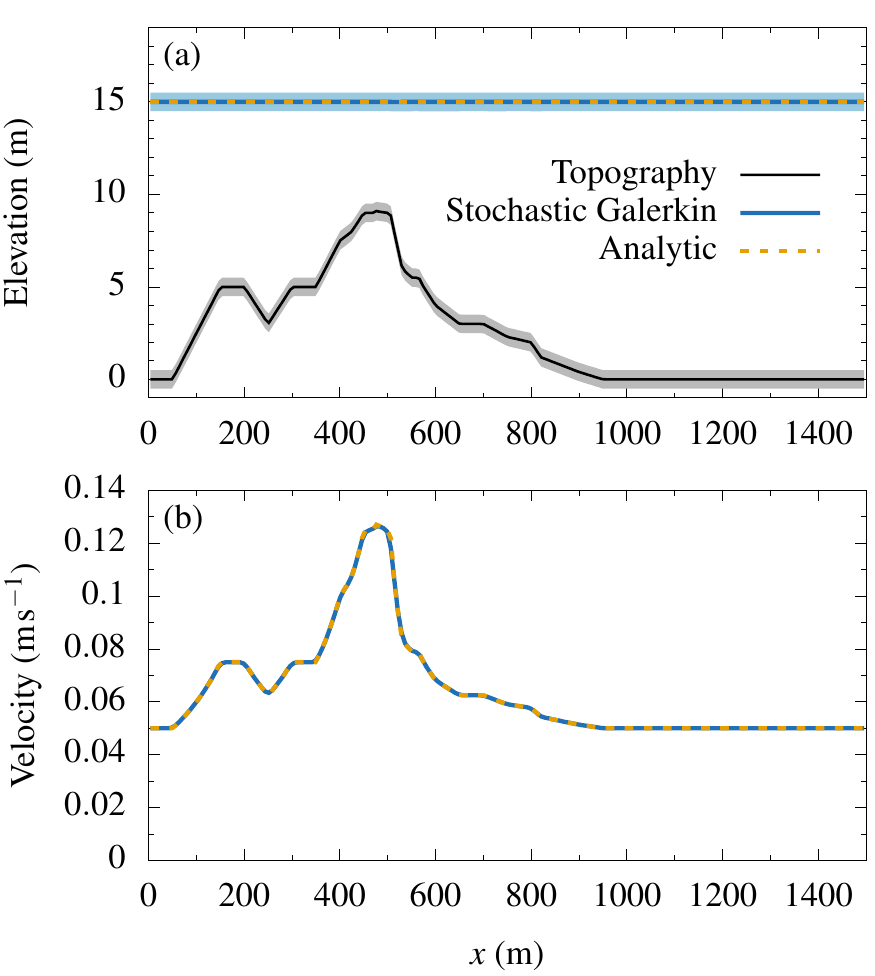}
    \end{subfigure}
    \caption{\rev{Steady state solution over an irregular and uncertain bed.
    (a) Bed elevation and free-surface elevation.
    Mean values are marked with solid lines, with shading indicating plus or minus one standard deviation of uncertainty.
    (b) Spatial profile of depth-average velocity.
    The mean velocity is marked with a solid line and the standard deviation of velocity is negligible.
    In both panels, the analytic mean solution is marked with a dashed yellow line.}}
    \label{fig:tsengSteadyState-flow}
\end{figure}

\rev{For this purpose, following \citet{tseng2004}, the 1D domain is $[\SI{0}{\meter}, \SI{1500}{\meter}]$ with $M = 200$ elements such that the mesh spacing is $\Delta x = \SI{7.5}{\meter}$.
The mean topography has an irregular profile as specified by \citet{goutal-maurel1997} and is shown in figure~\ref{fig:tsengSteadyState-flow:elevation}.
The standard deviation of topography is $\sigma_z(x) = \SI{0.5}{\meter}$ across the entire domain, which is chosen since it provides an acceptable margin of error for local-scale flood modelling applications \citep{schumann-bates2018}.
The initial free-surface elevation is \SI{15}{\meter} with no uncertainty.
Subcritical boundary conditions are imposed such that the mean upstream discharge per unit-width is \SI{0.75}{\meter\squared\per\second} and the mean downstream water depth is \SI{15}{\meter}, with no uncertainty on either boundary condition.
Transmissive boundary conditions are used for the upstream water depth and downstream discharge.}
\rev{The test uses a timestep of $\Delta t = \SI{0.5}{\second}$ and terminates at $t = \SI{100000}{\second}$ when the water has converged on a steady state down to a convergence error of \SI{e-8}{\meter} as measured by equation~\eqref{eqn:convergence}.}

\rev{The stochastic Galerkin model is configured with a basis order $P=3$, and the resulting steady-state solution is shown in figure~\ref{fig:tsengSteadyState-flow}.
Since the flow is very slow, the free-surface elevation is uniform.
The free-surface elevation has the same standard deviation as the topography, which is due to the water depth at the downstream boundary being fixed at \SI{15}{\meter} with no uncertainty.
As a result, the water depth profile is certain, and the uncertainty in the free-surface elevation is determined entirely by the topography.}
\rev{The corresponding depth-averaged velocity profile is shown in figure~\ref{fig:tsengSteadyState-flow:u}.
The stochastic Galerkin model calculates a velocity profile that is in excellent agreement with the analytic solution.
The predicted velocity profile is certain with negligible standard deviation.
This is expected since the water depth profile is unaffected by the prescribed topographic uncertainty.
These results verify the robustness of the stochastic Galerkin model for flow over a highly-irregular, uncertain bed.}

%% file: conclusions.tex
\section{Conclusions and outlook}

This paper proposed a stochastic Galerkin shallow flow model using the surface gradient method to ensure well-balancing.
The well-balancedness of the stochastic model was studied theoretically \rev{and verified in a series of idealised numerical tests, with a prescribed uncertainty that characterises topographic or bathymetric measurement errors.}
\rev{These tests demonstrated that the stochastic Galerkin model preserves discrete balance between flux and topography gradients to machine precision, and accurately reproduces steady flow over highly irregular and uncertain topography.}
\rev{Despite its relatively low-order Wiener-Hermite basis, the stochastic Galerkin model was able to simulate strongly nonlinear behaviour, producing discontinuous and highly non-Gaussian probability distributions that compared favourably against a Monte Carlo reference solution.}

It is likely that a more sophisticated treatment using a multiresolution wavelet decomposition \citep{lemaitre2004a,pettersson2014} or multi-element discretisation of the stochastic dimension \citep{wan-karniadakis2006,li-stinis2015} would further improve the representation of highly non-Gaussian and discontinuous probability distributions.
It is also likely that a decomposition or discretisation of stochastic space would simplify a stochastic reformulation of positivity preserving schemes involving nonlinear operators, so enabling the development of stochastic wetting-and-drying processes \rev{that are essential for flood inundation modelling}.
Such improvements are the subject of future work.

\rev{The one-dimensional stochastic Galerkin model presented here is ideally suited to river routing applications:
using as few as two expansion coefficients per element, the stochastic model can account for river bed measurement errors to produce a probabilistic flow simulation that is accurate and computationally efficient.}
Without parallelisation, the stochastic Galerkin model is \rev{at least} 100 times faster than a Monte Carlo simulation.
With parallelisation on commodity hardware, the computation time for the stochastic Galerkin model would be on par with a single iteration of the deterministic solver.

%The first test verified that a lake-at-rest over an uncertain bed was preserved to machine precision.
%The second test challenged the stochastic Galerkin model using steady-state critical flow over an uncertain bed.
%The test was designed to provoke strongly nonlinear behaviour such that about 50\% of flows were subcritical and the other 50\% develop a transcritical shock.
%The resulting probability distributions were discontinuous and highly non-Gaussian.
%\rev{A third test verified that the stochastic Galerkin model accurately simulates steady flow over a highly irregular and uncertain bed profile.}

%% file: ack.tex
\section*{Acknowledgements}

This work is part of the SEAMLESS-WAVE project (SoftwarE infrAstructure for Multi-purpose fLood modElling at variouS scaleS based on WAVElets) which is funded by the UK Engineering and Physical Sciences Research Council (EPSRC) grant EP/R007349/1.
For information about the SEAMLESS-WAVE project visit \url{https://www.seamlesswave.com}.
\rev{The authors are grateful to the editors and two anonymous reviewers for their helpful comments.}

%% file: software.tex
\section{Instructions for running the shallow flow models}

Python 3 software is available for download from Zenodo \citep{shaw-kesserwani2019b} which includes
\begin{enumerate}
    \item a one-dimensional shallow flow model which can operate as a stochastic Galerkin model or a deterministic model
    \item code for running Monte Carlo iterations of the deterministic model
    \item preconfigured lake-at-rest, steady-state critical flow \rev{and steady-state flow over an irregular bed} test cases
\end{enumerate}
First ensure that Python 3, NumPy and SciPy are installed.
Then, install the shallow flow model:
\begin{verbatim}
pip3 install --user --editable swepc.python
\end{verbatim}

\subsection{Running the Stochastic Galerkin Model}
To run the stochastic Galerkin model:
\begin{Verbatim}[commandchars=\\\{\}]
swepc --degree \textit{<P>} \textit{<testCase>} \textit{<discretisation>} --output-dir \textit{<directory>}
\end{Verbatim}
where \Verb[commandchars=\\\{\}]+\textit{<testCase>}+ is either \Verb+lakeAtRest+, \rev{\Verb+criticalSteadyState+ or \Verb+tsengSteadyState+},
\Verb[commandchars=\\\{\}]+\textit{<discretisation>}+ is either \Verb+wellBalancedH+ for the well-balanced surface gradient method, or \Verb+centredDifferenceH+ for the centred difference method given by equations~\eqref{eqn:flux-centred} and \eqref{eqn:source-centred}.
If \Verb[commandchars=\\\{\}]+--degree \textit{<P>}+ is omitted then the default basis order $P = 3$ is used.
If \Verb[commandchars=\\\{\}]+--degree 0+ is specified then the model operates as a deterministic model.
On running the model, the following text files are written to the specified output directory:
\begin{description}
\item[\Verb+coefficients.dat+]{Stochastic Galerkin expansion coefficients for topography $z$, water \rev{depth} $h$ and discharge $q$}
\item[\Verb+statistics.dat+]{Mean, standard deviation, skew and kurtosis calculated using equations~\eqref{eqn:moment}, \eqref{eqn:mean} and \eqref{eqn:variance}}
\item[\Verb+derived-statistics.dat+]{Statistics for the \rev{depth-averaged velocity} and the free-surface elevation, a derived quantity calculated from equation~\eqref{eqn:h-eta-z}}
\end{description}
Each file has one line per element with explanatory header rows prefixed by \Verb+#+.

\subsection{Calculating Probability Density Functions}
A probability density function for any variable can be calculated from the expansion coefficients for a given element:
\begin{Verbatim}[commandchars=\\\{\}]
sed -n \textit{<line>}p coefficients.dat
   | swepdf --min \textit{<value>} --max \textit{<value>} \textit{<variable>} > pdf.dat
\end{Verbatim}
where \Verb[commandchars=\\\{\}]+\textit{<line>}+ is the line number in \Verb+coefficients.dat+ corresponding to the chosen element, and \Verb[commandchars=\\\{\}]+\textit{<variable>}+ is \Verb+z+ (topography), \Verb+water+ (water \rev{depth}), \Verb+q+ (discharge) or \Verb+derived-eta+ (free-surface elevation).
The probability density function is calculated between the \Verb+--min+ and \Verb+--max+ values and written to \Verb+pdf.dat+. 

\subsection{Running Monte Carlo Simulations}
To run a Monte Carlo simulation:
\begin{Verbatim}[commandchars=\\\{\}]
swepc --monte-carlo --mc-iterations \textit{<value>}
  criticalSteadyState \textit{<discretisation>} --output-dir \textit{<directory>}
\end{Verbatim}
which writes the following text files:
\begin{description}
\item[\Verb+statistics.dat+]{Mean, standard deviation, skew and kurtosis for topography, water \rev{depth} and discharge}
\item[\Verb+derived-statistics.dat+]{Mean and standard deviation for free-surface elevation \rev{and depth-averaged velocity}}
\cprotect\item[\Verb[commandchars=\\\{\}]+sample\textit{<n>}.dat+]{Deterministic model data with one file per element (numbered $0, \ldots, M-1$), each having one line per Monte Carlo iteration}
\end{description}